\documentclass[10pt]{amsart}

\usepackage{hyperref}

\hypersetup{ colorlinks=true, linkcolor=blue, filecolor=magenta, urlcolor=cyan, }

\usepackage{graphics}
\usepackage{graphicx}
\usepackage{epsfig}
\usepackage{amsmath}
\usepackage{graphics}
\usepackage{graphicx}
\usepackage{epstopdf}
\usepackage{multirow}
\usepackage{booktabs}
\usepackage{ifpdf}
\usepackage{color}
\usepackage{amssymb}
\usepackage{amsthm}
\usepackage{amsmath}
\usepackage{matlab-prettifier}
\usepackage{tikz-cd,cleveref}
\usepackage{todonotes}
\usepackage{nicematrix}
\usepackage{ytableau}
\usepackage{tikz}
\usepackage{graphicx}
\usepackage{comment}
\usepackage[T1]{fontenc}
\usepackage{adjustbox}
\usepackage{fancybox}
\usetikzlibrary{decorations.pathreplacing,
	calligraphy,% must be loaded after decorations.pathreplacing
	matrix}
\usepackage{float}
\restylefloat{table}

\theoremstyle{plain}% Theorem-like structures provided by amsthm.sty
\newtheorem{thm}{Theorem}[section]
\newtheorem{pro}[thm]{Proposition}

\theoremstyle{definition}
\newtheorem{dfn}[thm]{Definition}

\newtheorem{rem}[thm]{Remark}

\newtheorem{question}[thm]{Question}
\usepackage{lineno}
\numberwithin{equation}{section}

\cornersize{0.2}
\newcommand\mlnode[1]{\ovalbox{\begin{tabular}{@{}c@{}}#1\end{tabular}}}

\def \p {\mathbb{P}}

\def \z {\mathbb{Z}}

\def \c {\mathbb{C}}

\def \b {\mathbb{B}}

\def \codim {\operatorname{codim}}

\def \det {\operatorname{det}}
\def \rk {\operatorname{\mathbf{R}}}
\def \brk {\operatorname{\underline{\mathbf{R}}}}

\def \sgn {\operatorname{sgn}}

\def \SL {\operatorname{SL}}
\def \sym {\operatorname{Sym}}

\begin{document}

\bibliographystyle{plain}

\title{The Border Rank of the $4 \times 4$ Determinant Tensor is Twelve}

\author{Jong In Han}
\address{%Jong In Han, 
School of Mathematics, Korea Institute for Advanced Study (KIAS), 85 Hoegi-ro, Dongdaemun-gu, Seoul 02455, Republic of Korea}
\email{jihan09@kias.re.kr}

\author{Jeong-Hoon Ju}
\address{%Jeong-Hoon Ju, 
Department of Mathematical Sciences, University of Copenhagen, 5 Universitetsparken, 2100 Copenhagen Ø, Denmark}
\email{jeonghoon.ju@math.ku.dk}

\author{Yeongrak Kim}
\address{%Yeongrak Kim, 
Department of Mathematics \& Institute of Mathematical Sciences, Pusan National University, 
	2 Busandaehak-ro 63beon-gil, Geumjeung-gu, 46241 Busan, Republic of Korea}
\email{yeongrak.kim@pusan.ac.kr}

\thanks{}

\begin{abstract}
	We show that the border rank of the $4 \times 4$ determinant tensor is at least $12$ over $\c$, using the fixed ideal theorem introduced by  Buczy\'nska--Buczy\'nski and the method by Conner--Harper--Landsberg. Together with the known upper bound, this implies that the border rank is exactly $12$.
\end{abstract}



\maketitle
%

%% main text

% Introduction
\section{Introduction}
Let $V_1,\cdots, V_n$ be finite dimensional $\c$-vector spaces. A nonzero tensor $T \in V_1 \otimes \cdots \otimes V_n$ is of \emph{rank one} if $T = v_1 \otimes \cdots \otimes v_n$ for some $v_i \in V_i$. The \emph{tensor rank} of $T$, denoted by $\mathbf{R} (T)$, is the smallest integer $r$ such that $T$ can be written as a sum of $r$ rank one tensors. It generalizes the usual notion of the matrix rank, however, it does not satisfy a certain semi-continuity, and hence it is not comfortable to apply geometric ideas to study the tensor rank. A more geometric notion would be the \emph{border rank} of $T$, denoted by $\underline{\mathbf{R}} (T)$, which is the smallest integer $r$ such that $T$ can be written as a limit of a sum of $r$ rank one tensors. Together with the $r$-th secant variety of the Segre variety of rank one tensors, we may deal with the border rank in an algebro-geometric way. 

The main object in this paper is the $4 \times 4$ determinant tensor $\det_4 \in \c^4 \otimes \c^4 \otimes \c^4 \otimes \c^4$ considered it as a multilinear map. Note that many important problems in complexity theory concern the determinant polynomial of a square matrix of variables, however, this tensor is also a fundamental object in algebra and geometry, and have a better prospect to understand its rank complexities. Note that the usual determinant polynomial $\operatorname{Det}_n$, considered as a symmetric tensor, is the same as the Kronecker square $\det_n^{\boxtimes 2}$. Studying the border rank of a tensor and its Kronecker powers is a very interesting problem, particularly due to a connection to the exponent of matrix multiplication together with Strassen's laser method \cite{strassen1987laser}. We refer to \cite{conner2022kronecker} for more detailed discussions, and we expect that the study of the tensor and border ranks of $\det_n$ and its Kronecker powers would be a foundational work in geometric complexity theory with these reasons.

Let us briefly summarize what are known about the tensor rank and the border rank of $\det_4$. 
It was already known that $\mathbf{R} (\det_4) \le 12$ by addressing an explicit decomposition as a sum of $12$ rank one tensors \cite{MR4847254}. In the previous paper of the authors, it was shown that $\mathbf{R} (\det_4) = 12$ and $11 \le \underline{\mathbf{R}} (\det_4)$, by using a method named the recursive Koszul flattening method \cite{MR3987862, RKF}. Hence, the only remaining problem is to determine whether the border rank $\underline{\mathbf{R}}(\det_4)$ is $11$ or $12$.

Apolarity theory is very helpful when we want to analyze rank of a specific tensor or a polynomial. Classical apolarity theory, extending a work of Sylvester, associates the ideal of derivations (called the \emph{apolar ideal} of $F$) which annihilates a given homogeneous polynomial $F$, and it indicates a way how to find a decomposition of $F$ as a sum of powers of linear forms. This theory can be naturally generalized in multigraded setting \cite{galazka2023apolarity}. Note that the analogous apolarity theory for border rank is also introduced by Buczy\'nska--Buczy\'nski \cite{MR4332674}. In particular, if a given tensor $T$ can be written as a limit of a sum of $r$ rank one tensors (i.e., $\underline{\mathbf{R}}(T) \le r$), the border apolarity theorem \cite[Theorem 3.15]{MR4332674} assures the existence of a (multi)homogeneous ideal $I$ contained in the apolar ideal of $T$ which can be obtained as a scheme of limits of ideals of $r$ points (slip). Furthermore, if there is a group $G$ acting on the Segre variety and preserving $T$, one can find such $I$ that is invariant under the action of a Borel subgroup of $G$ \cite[Theorem 4.3]{MR4332674}. This fixed ideal theorem generalizes the normal form lemma of Landsberg--Micha{\l}ek \cite[Lemma 3.1]{michalek2017normalform}, and shows how effective the algebro-geometric approach is in studying border ranks. Using this theory, Conner--Harper--Landsberg described an algorithm to enumerate a set of parametrized families of ideals which could satisfy the conclusion of the fixed ideal theorem \cite{MR4595287}. If this enumeration fails, then the assumption $\underline{\mathbf{R}}(T) \le r$ also fails, so this method might be helpful to improve a lower bound on the border rank of a tensor with large symmetries. 

By using these methods, we determine the border rank of $\det_4$:
\begin{thm}\label{ThmMain}
	The border rank of $\det_4$ is $12$ over any subfield of $\c$.
\end{thm}

The structure of the paper is as follows. In Section \ref{Sect:Preliminaries}, we review some basic terminology on tensors and border apolarity theory. And then we describe an explicit border rank criterion for concise tensors of order $4$, which is an analogue of the algorithm provided in \cite{MR4595287}. In Section \ref{Sect:Proof}, we perform the test for $\det_4$ and conclude that its border rank cannot be $11$. Also, we give some remarks on the fixed ideal theorem and a question. Some of the calculations are computer-assisted, and the code is available at \cite{hanjukim2025appendix}.

%--------------------------------------------------------------

%--------------------------------------------------------------
\section*{Acknowledgements}
J. I. H. was partially supported by the National Research Foundation of Korea (NRF) grant funded by the Korea government (MSIT) (No. RS-2024-00352592) when he was at Korea Advanced Institute of Science and Technology (KAIST). J. I. H. is currently supported by a KIAS Individual Grant (MG101401) at Korea Institute for Advanced Study.
J.-H. J. and Y. K. are supported by the Basic Science Program of the NRF of Korea (NRF-2022R1C1C1010052) and by the Basic Research Laboratory (grant MSIT no. RS-202400414849). J.-H. J. is currently supported by the Novo Nordisk Foundation (grant NNF20OC0059939).
The authors thank J. M. Landsberg for helpful comments and discussions.

%--------------------------------------------------------------

\section{Preliminaries}\label{Sect:Preliminaries}

Throughtout this paper, the base field is the complex number field $\c$ if not stated otherwise, and we use the following notations:
\begin{itemize}
	\item $V, V_i$ : finite dimensional vector spaces over $\c$;
	\item $V^*$ : the dual vector space of $V$;
	\item $S^dV$ : the space of symmetric tensors of order $d$ in $V^{\otimes d}$, or equivalently the space of homogeneous polynomials on $V$ of degree $d$ and zero;
	\item $\sym(V):=\oplus_{d\geq 0}S^dV$;
	\item $\Lambda^dV$ : the space of skew-symmetric tensors of order $d$ in $V^{\otimes d}$;
	\item For a Young tableau $\lambda$, $S_{\lambda}V$ denotes the Schur module corresponding to $\lambda$;
	\item $\c^*:=\c \setminus \{0\}$;
	\item $[d] = \{1, 2, \cdots, d \}$ where $d$ is a positive integer;
	\item $\mathfrak{S}_d$ : the symmetric group on $d$ letters;
	\item $\operatorname{SL_n}(\text{resp.}~\mathfrak{sl}_n)$ : special linear group (resp. algebra) of $\c^n$.
\end{itemize}

\subsection{Concise Tensor}\label{SubSect:Tensor}

An element in $V_1 \otimes V_2 \otimes \cdots \otimes V_d$ is called a \emph{tensor} of \emph{order} $d$. What we consider in this paper mainly is the \textit{$4 \times 4$ determinant tensor} defined as
\begin{equation}
	\det_4=\sum_{\sigma \in \mathfrak{S}_4}\sgn(\sigma) e_{\sigma(1)} \otimes e_{\sigma(2)} \otimes e_{\sigma(3)} \otimes e_{\sigma(4)} \in V_1 \otimes V_2 \otimes V_3 \otimes V_4
\end{equation} 
where $V_1 \cong V_2 \cong V_3 \cong V_4 \cong \c^4$, which is of order $4$. 

For a given tensor $T \in V_1 \otimes V_2 \otimes \cdots \otimes V_d$ for $d \geq 2$, we can consider it as a linear map $$T_{V_1}:V_1^* \rightarrow V_2 \otimes \cdots \otimes V_d.$$ Let $T(V_1^*)$ denote the image $T_{V_1}(V_1^*)$ in $V_2 \otimes \cdots \otimes V_d$. Similarly, we can consider $T$ as a linear map $$T_{V_i}:V_i^* \rightarrow V_1 \otimes \cdots \otimes \hat{V_i} \otimes \cdots \otimes V_d,$$ and let $T(V_i^*)$ be its image.

\begin{dfn}[Concise tensor]
	A tensor $T \in V_1 \otimes \cdots \otimes V_d$ is said to be \emph{concise} if all the $T_{V_i}$'s are injective.
\end{dfn}

In other words, a tensor $T  \in V_1 \otimes \cdots \otimes V_d$ is concise if and only if for any $i$ there is no proper subspace $W_i$ of $V_i$ such that $$T \in V_1 \otimes \cdots \otimes V_{i-1} \otimes W_i \otimes V_{i+1} \otimes \cdots \otimes V_d.$$ One can easily check that the main object $\det_4$ is concise.

\subsection{Border Apolarity Lemma and Fixed Ideal Theorem}\label{SubSect:Border Apolarity Lemma and Fixed Ideal Theorem}

In this section, we review the definition of border rank, and then introduce some important results in border apolarity theory.

\begin{dfn}[Secant variety and border $X$-rank]
	Let $X \subset \p V$ be a nondegenerate projective variety. The \emph{$r$-th secant variety of $X$} is
	\begin{equation*}
		\sigma_r(X)=\overline{\bigcup_{p_1,\cdots,p_r \in X} \langle p_1,\cdots,p_r \rangle}.
	\end{equation*}
	For a point $p$ in $\p V$, its \emph{border $X$-rank}, denoted by $\brk_X(p)$, is defined as
	\begin{equation*}
		\brk_X(p)=\min\{r \in \z~|~p \in \sigma_r(X)\}.
	\end{equation*}
\end{dfn}

\begin{dfn}[Border rank]
	If $X$ is a Segre variety  $$\operatorname{Seg}(\p V_1 \times \p V_2 \times \cdots \times  \p V_d) \subset \p(V_1 \otimes V_2 \otimes \cdots \otimes V_d)$$ and $T \in V_1 \otimes V_2 \otimes \cdots \otimes V_d$, then the border $X$-rank of $[T] \in  \p(V_1 \otimes V_2 \otimes \cdots \otimes V_d)$ is called the \emph{border rank} of $T$, and simply denoted by  $\brk(T)$.
\end{dfn}

One of the useful tools to investigate the border rank is \emph{(border) apolarity theory}, which compares the annihilator and defininig ideal of a given tensor. An element in $\sym(V^*)$ can be considered as a differential operator on elements in $V$. Hence, an element $$\Theta \in \sym(V_1^*) \otimes \cdots \otimes  \sym(V_d^*)$$ also can be considered as a differential operator on elements $T \in V_1 \otimes \cdots \otimes V_d$. We denote the differentiation of $\Theta$ to $T$ as $\Theta \mathbin \lrcorner T$.

\begin{dfn}[Annihilator]\label{DefAnn}
	An annihilator of a tensor $T \in V_1 \otimes \cdots \otimes V_d$, denoted by $\operatorname{Ann}(T)$, is defined by
	\begin{equation*}
		\operatorname{Ann}(T)=\{\Theta \in \sym(V_1^*) \otimes \cdots \otimes  \sym(V_d^*)~|~\Theta \mathbin \lrcorner T=0\}.
	\end{equation*}
\end{dfn}

\begin{rem}\label{RmkPerp}
	In Definition \ref{DefAnn}, all elements $\Theta$ of multidegree $(i_1,\cdots,i_d)$ with one of $i_1,\cdots,i_d$ greater than $1$ are all in $\operatorname{Ann}(T)$.  We describe elements $\Theta$ of multidegree $(i_1,\cdots,i_d)$ with $i_1,\cdots,i_d$ all less than or equal to $1$ in terms of linear annihilator and flattening. For a subset $X \subset V$, its \emph{linear annihilator}, denoted by $X^{\perp}$, is defined as $$X^{\perp}=\{f \in V^*~|~f(x)=0~\text{for all}~x \in X\}.$$ As we could consider $T \in V_1 \otimes \cdots \otimes V_d$ as a linear map $T_{V_i}$, we also are able to consider $T$ as a linear map $T_{V_{j_1}\otimes \cdots \otimes V_{j_s}}$ called \emph{flattening} of $T$ and its image $T(V_{j_1}^* \otimes \cdots \otimes  V_{j_s}^*)$ for $s \in [d]$ and pairwise distinct $j_1,\cdots,j_s \in [d]$.  The elements of $\Theta$ of multidegree $(i_1,\cdots,i_d)$ with $i_k\le 1$ for all $k$ are exactly the elements in one of $T^{\perp}$ or $T(V_{j_1}^* \otimes \cdots \otimes  V_{j_s}^*)^{\perp}$'s for  $s \in [d]$ and  pairwise distinct $j_1,\cdots,j_s \in [d]$. Note that $T(V_{j_1}^* \otimes \cdots \otimes V_{j_{d-1}}^*)^{\perp}$ is the zero space when $T$ is concise.
\end{rem}

\begin{rem}\label{RmkConciseIdeal}
	If $T$ is concise, then $T(V_{j_1}^* \otimes \cdots \otimes V_{j_{d-1}}^*)^{\perp}$ is zero, and so any $I$ such that $I \subset \operatorname{Ann}(T)$ must satisfy $$I_{10\cdots 0}=I_{010\cdots 0}=\cdots=I_{0\cdots 01}=0.$$
\end{rem}

We state two main theorems on border apolarity theory by using multihomogeneous ideal and symmetry group of a tensor.  An ideal $I \subset \sym(V_1^*) \otimes  \cdots \otimes \sym(V_d^*)$ is said to be \textit{multihomogeneous} if it is homogeneous on each factor $\sym(V_i^*)$. Each element in the multihomogeneous ideal $I$ has multidegree $(i_1,\cdots,i_d)$ for some $(i_1,\cdots,i_d)\in (\z_{\geq 0})^{\times d}$, and we let $I_{i_1,\cdots, i_d}$ denote the degree $(i_1,\cdots,i_d)$ part of $I$, which is a subspace of $S^{i_1}V_1^* \otimes \cdots \otimes S^{i_d}V_d^*$. For $T \in V_1 \otimes \cdots \otimes V_d$, define its \emph{symmetry group} as
\begin{equation*}
	G_T:=\{(g_{V_1},\cdots,g_{V_d}) \in \operatorname{GL}(V_1) \times \cdots \times  \operatorname{GL}(V_d)/(\c^*)^{\times (d-1)}~|~(g_{V_1},\cdots,g_{V_d}) \cdot T=T\}.
\end{equation*}
Here, we took the quotient by $$(\c^*)^{\times (d-1)}=\{(\lambda_1 \operatorname{Id}_{V_1},\cdots,\lambda_d \operatorname{Id}_{V_d})~|~\lambda_1\cdots\lambda_d=1\}$$ because $$\{(\lambda_1 \operatorname{Id}_{V_1},\cdots, \lambda_{d-1} \operatorname{Id}_{V_{d-1}},\frac{1}{\lambda_1\cdots\lambda_{d-1}} \operatorname{Id}_{V_d})~|~(\lambda_1,\cdots,\lambda_{d-1}) \in (\c^*)^{\times (d-1)}\}$$ is the kernel of the map $$\operatorname{GL}(V_1) \times \cdots \times \operatorname{GL}(V_d) \rightarrow \operatorname{GL}(V_1 \otimes \cdots \otimes V_d).$$ 

Suppose $T$ has a border rank $r$ decomposition $T=\lim_{t\to 0}\sum^r_{i=1}T_i(t)$. Then for each $t\neq 0$, we get a $\z^d$-graded ideal $I_t\subset \sym(V_1^*)\otimes \cdots\otimes \sym(V_d^*)$ corresponding to the points $T_1(t), \cdots, T_r(t)$.
Assume that
\begin{equation}\label{codim_condition}
\begin{split}
	\text{$I_{t,i_1,\cdots,i_d}\subset S^{i_1}V_1^*\otimes \cdots\otimes S^{i_d}V_d^*$ has codimension $r$}\\
	\text{whenever $\dim (S^{i_1}V_1^*\otimes \cdots\otimes S^{i_d}V_d^*)\geq r$.}\hspace{0.6cm}
\end{split}
\end{equation}

In this case, $I_{t,i_1,\cdots,i_d}$ can be considered as an element in the Grassmannian $$G(\dim(S^{i_1}V_1^*\otimes \cdots\otimes S^{i_d}V_d^*)-r,S^{i_1}V_1^*\otimes \cdots\otimes S^{i_d}V_d^*).$$
Then we get an ideal $I\subset \sym(V_1^*)\otimes \cdots\otimes \sym(V_d^*)$ whose degree $(i_1,\cdots,i_d)$ part is $$I_{i_1,\cdots,i_d}=\lim_{t\to 0}I_{t,i_1,\cdots,i_d}$$ for all $(i_1,\cdots,i_d)$ such that $\dim (S^{i_1}V_1^*\otimes \cdots\otimes S^{i_d}V_d^*)\geq r$. Here, the limit is taken in $$G(\dim(S^{i_1}V_1^*\otimes \cdots\otimes S^{i_d}V_d^*)-r,S^{i_1}V_1^*\otimes \cdots\otimes S^{i_d}V_d^*).$$ For those $(i_1,\cdots,i_d)$ such that $\dim (S^{i_1}V_1^*\otimes \cdots\otimes S^{i_d}V_d^*)< r$, we let $I_{i_1,\cdots,i_d}=0$.
We call $I$ an \textit{ideal corresponding to a border rank $r$ decomposition of $T$}. Note that we define an ideal corresponding to a border rank decomposition only for those satisfying \eqref{codim_condition}.
\begin{rem}
	The condition \eqref{codim_condition} does not always need to be satisfied even if $T_1(t)+\cdots+T_r(t)$ are minimal decompositions for each $t$. For example, consider the following decomposition of $\det_4\in V_1\otimes V_2\otimes V_3\otimes V_4$ given in \cite{MR4847254}.
    This is a minimal decomposition by \cite{RKF}.
	\begin{align*}
		\det_4 =\frac{1}{2}(
		&(e_1-e_2)\otimes(e_3-e_4)\otimes(e_3+e_4)\otimes(e_1+e_2)\\
		&-(e_1-e_3)\otimes(e_2-e_4)\otimes(e_2+e_4)\otimes(e_1+e_3)\\
		&+(e_1-e_4)\otimes(e_2-e_3)\otimes(e_2+e_3)\otimes(e_1+e_4)\\
		&+(e_2-e_3)\otimes(e_1-e_4)\otimes(e_1+e_4)\otimes(e_2+e_3)\\
		&-(e_2-e_4)\otimes(e_1-e_3)\otimes(e_1+e_3)\otimes(e_2+e_4)\\
		&+(e_3-e_4)\otimes(e_1-e_2)\otimes(e_1+e_2)\otimes(e_3+e_4)\\
		&+(e_1+e_2)\otimes(e_3+e_4)\otimes(e_3-e_4)\otimes(e_1-e_2)\\
		&-(e_1+e_3)\otimes(e_2+e_4)\otimes(e_2-e_4)\otimes(e_1-e_3)\\
		&+(e_1+e_4)\otimes(e_2+e_3)\otimes(e_2-e_3)\otimes(e_1-e_4)\\
		&+(e_2+e_3)\otimes(e_1+e_4)\otimes(e_1-e_4)\otimes(e_2-e_3)\\
		&-(e_2+e_4)\otimes(e_1+e_3)\otimes(e_1-e_3)\otimes(e_2-e_4)\\
		&+(e_3+e_4)\otimes(e_1+e_2)\otimes(e_1-e_2)\otimes(e_3-e_4))
	\end{align*}
	Consider this as a border rank 12 decomposition, i.e., $T_i(t)\in V_1\otimes V_2\otimes V_3\otimes V_4$ is the point corresponding to the $i$-th row in the above equation for all $t$.
	Then $I_{t,1100}\subset V_1^*\otimes V_2^*$ does not have codimension 12 since the projections of $\{T_i(t)\}_{i=1,\cdots,12}$ to $V_1\otimes V_2$ are not linearly independent.
	Indeed, as the span of
	\begin{align*}
		&(e_1-e_2)\otimes(e_3-e_4),
		(e_1-e_3)\otimes(e_2-e_4),\\
		&(e_1-e_4)\otimes(e_2-e_3),
		(e_2-e_3)\otimes(e_1-e_4),\\
		&(e_2-e_4)\otimes(e_1-e_3),
		(e_3-e_4)\otimes(e_1-e_2),\\
		&(e_1+e_2)\otimes(e_3+e_4),
		(e_1+e_3)\otimes(e_2+e_4),\\
		&(e_1+e_4)\otimes(e_2+e_3),
		(e_2+e_3)\otimes(e_1+e_4),\\
		&(e_2+e_4)\otimes(e_1+e_3),
		(e_3+e_4)\otimes(e_1+e_2)
	\end{align*}	
	has dimension $9$ in $V_1\otimes V_2$, the codimension of $I_{t,1100}\subset V_1^*\otimes V_2^*$ is $9$.
    Still, the fixed ideal theorem (\Cref{thmFixedIdeal}) guarantees the existence of a border rank $r$ decomposition satisfying \eqref{codim_condition} when $\brk(T)\le r$.
\end{rem}

Buczy\'nska--Buczy\'nski showed the weak border apolarity theorem holds:

\begin{thm}[Weak border apolarity theorem, \cite{MR4332674}]\label{thmWeakApolar}
	Let $T \in V_1 \otimes \cdots \otimes V_d$. If $\brk(T) \leq r$, then there exists a multihomogeneous ideal $$I \subset \sym(V_1^*) \otimes  \cdots \otimes \sym(V_d^*)$$ satisfying the following:
	\begin{itemize}
		\item [(\romannumeral1)] $I \subset \operatorname{Ann}(T)$;
		\item [(\romannumeral2)] For each $(i_1,\cdots,i_d) \in (\z_{\geq 0})^{\times d}$, the codimension of $I_{i_1,\cdots, i_d}$ as a subspace in $S^{i_1}V_1^* \otimes \cdots \otimes S^{i_d}{V_d^*}$ is
		\begin{equation}\label{eqCodim}
			\operatorname{codim}(I_{i_1,\cdots, i_d}) = \min(r,\dim (S^{i_1}V_1^* \otimes \cdots \otimes S^{i_d}V_d^*)).
		\end{equation}
	\end{itemize}
\end{thm}

Furthermore, they showed even a stronger statement holds: for a conneted solvable group $H\subset G_T$, there exists an $H$-fixed ideal corresponding to a border rank $r$ decomposition of $T$.

\begin{thm}[Fixed ideal theorem, \cite{MR4332674}]\label{thmFixedIdeal}
	Let $T \in V_1 \otimes \cdots\otimes V_d$, and let $H \subset G_T$ be a connected solvable group. If $\brk(T) \leq r$, then there exists an ideal $$I \subset \sym(V_1^*) \otimes\cdots \otimes \sym(V_d^*)$$ corresponding to a border rank $r$ decomposition of $T$ that is $H$-fixed, i.e., $$b \cdot I_{i_1,\cdots, i_d}=I_{i_1,\cdots, i_d}$$ for all $b \in H$ and all multidegree $(i_1,\cdots,i_d)$.
\end{thm}

\subsection{Border Rank Criterion}\label{SubSect:IBorderRankCriterion}

In this section, we see how one can check whether a tensor may have border rank at most $r$. This method is introduced by Conner--Harper--Landsberg  \cite{MR4595287} using the fixed ideal theorem (\Cref{thmFixedIdeal}) systemically. This is the state-of-art method to give a lower bound of border ranks, and they induced new lower bounds for various tensors of order $3$. 
Their method works for higher order tensors similarly.
We state the method for order 4 tensors to deal with the $4 \times 4$ determinant tensor $\det_4$.

Let $T \in V_1 \otimes V_2 \otimes V_3 \otimes V_4$ be a concise tensor with  $\dim V_1 \leq \dim V_2 \leq \dim V_3 \leq \dim V_4$ and $\mathbb{B}\subset G_T$ be a connected solvable group, for instance a Borel subgroup.
Assume that $\brk(T)\leq r$ so that there exists a border rank $r$ decomposition, and let $I$ be the multihomogeneous ideal corresponding to the border rank $r$ decomposition of $T$ that is $\b$-fixed obtained in \Cref{thmFixedIdeal}. Then the following must hold:
\begin{itemize}
	\item [(\romannumeral1)] $I \subset \operatorname{Ann}(T)$ holds. That is, by Remark \ref{RmkPerp},
	\begin{itemize}
		\item [$\bullet$] $I_{1100} \subset T(V_3^* \otimes V_4^*)^{\perp},I_{1010} \subset T(V_2^* \otimes V_4^*)^{\perp},\cdots,I_{0011} \subset T(V_1^* \otimes V_2^*)^{\perp}$;
		\item [$\bullet$] $I_{1110} \subset T(V_4^*)^{\perp},I_{1101} \subset T(V_3^*)^{\perp},I_{1011} \subset T(V_2^*)^{\perp},I_{0111} \subset T(V_1^*)^{\perp}$;
		\item [$\bullet$] $I_{1111} \subset T^{\perp} \subset V_1^* \otimes V_2^* \otimes V_3^* \otimes V_4^*$.
	\end{itemize}
	\item [(\romannumeral2)] For all $(i_1,i_2,i_3,i_4)$ such that $r \leq \dim (S^{i_1}V_1 ^*\otimes S^{i_2}V_2^* \otimes S^{i_3}V_3^* \otimes S^{i_4}V_4^*)$,  $\codim (I_{i_1,i_2,i_3,i_4}) = r$.
	\item [(\romannumeral3)] The multiplication map
	{\small
	\begin{equation*}
		\begin{aligned}
			(I_{(i_1-1),i_2,i_3,i_4} \otimes V_1^*) \oplus & (I_{i_1,(i_2-1),i_3,i_4} \otimes V_2^*) \oplus  (I_{i_1,i_2,(i_3-1),i_4} \otimes V_3^*) \oplus (I_{i_1,i_2,i_3,(i_4-1)} \otimes V_4^*)\\
			& \rightarrow  S^{i_1}V_1^*\otimes S^{i_2}V_2^* \otimes S^{i_3}V_3^* \otimes S^{i_4}V_4^*
		\end{aligned}		
	\end{equation*}
	}
	has its image in $I_{i_1,i_2,i_3,i_4}$ since $I$ is an ideal. 
	\item [(\romannumeral4)] Each $I_{i_1,i_2,i_3,i_4}$ is $\mathbb{B}$-fixed.
\end{itemize}
Therefore, if we show that there is no ideal $I$ satisfying (\romannumeral1)-(\romannumeral4) above, then we obtain $\brk(T)>r$.

\begin{rem}\label{RmkDetCodim}
	Let $V_1 \cong V_2 \cong V_3 \cong V_4 \cong \c^4$. For all $(i_1,i_2,i_3,i_4)$ with $i_1,\cdots,i_4$ at least two of them are greater than or equal to $1$, we obtain
	\begin{equation*}
		11 \leq \dim (S^{i_1}V_1 ^*\otimes S^{i_2}V_2^* \otimes S^{i_3}V_3^* \otimes S^{i_4}V_4^*).
	\end{equation*}
\end{rem}

We proceed to test whether there exists such an ideal $I$ or not as follows. At first, since $T$ is concise, let $$I_{1000}=I_{0100}=I_{0010}=I_{0001}=0$$ as in Remark \ref{RmkConciseIdeal}. Then take a $\mathbb{B}$-fixed subspace $F_{1100}$ of codimension $r$ in $V_1^* \otimes V_2^*$. If the following multiplication maps 
\begin{equation}\label{eq(2100)test}
	F_{1100} \otimes V_1^* \rightarrow S^2V_1^* \otimes V_2^*
\end{equation}
and
\begin{equation}\label{eq(1200)test}
	F_{1100} \otimes V_2^* \rightarrow V_1^* \otimes S^2V_2^*
\end{equation}
have images of codimension strictly less than $r$,
%, i.e., of dimension strictly greater than $\dim(V_1^* \otimes V_2^*)-r$,
then $F_{1100}$ cannot be a candidate of $I_{1100}$. These tests for maps (\ref{eq(2100)test}) and (\ref{eq(1200)test}) are respectively called $(2100)$-test and $(1200)$-test.  If for all $F_{1100}$, the multiplication maps have images of codimension strictly less than $r$, then there is no candidate of $I_{1100}$ and so we can conclude that $\brk(T)>r$. If there is a candidate of $I_{1100}$, we similarly take $F_{1010}$ and so on. Moreover, we can do $(2110),(1210),(1120)$-tests by taking a $\mathbb{B}$-fixed subspace $F_{1110}$ of codimension $r$ in $V_1^* \otimes V_2^* \otimes V_3^*$, and so on.

Suppose that we have some $4$-tuples $\{F_{1110},F_{1101},F_{1011},F_{0111}\}$ of candidates which pass $(2110),(1210),(1120)$-tests, $\cdots$, $(0211),(0121),(0112)$-tests, respectively. If for all $4$-tuples $\{F_{1110},F_{1101},F_{1011},F_{0111}\}$ of them have images of the map
\begin{equation}\label{eq(1111)test}
	(F_{1110} \otimes V_4^*) \oplus (F_{1101} \otimes V_3^*)\oplus(F_{1011} \otimes V_2^*)\oplus(F_{0111} \otimes V_1^*) \rightarrow V_1^* \otimes V_2^* \otimes V_3^* \otimes V_4^*
\end{equation}
of codimension strictly less than $r$, i.e., of dimension strictly greater than $\dim (V_1^* \otimes V_2^* \otimes V_3^* \otimes V_4^*)-r$, then there is no candidate of $I_{1111}$ and so we can conclude that $\brk(T) > r$. The test for the map (\ref{eq(1111)test}) is called the $(1111)$-test. We can also do $(1110)$-test for triples $\{F_{1100},F_{1010},F_{0110}\}$, but we explained $(1111)$-test as a representative because it is critical for $\det_4$. Moreover, if there is a $4$-tuple which passes the $(1111)$-test, then we can test similarly for higher order, for instance $(2111)$-test and $(2222)$-test. 

When the dimension of $F_{1110}$ is much bigger than that of $F_{1110}^{\perp}$, it is more convenient to perform the tests on the dual side using the following proposition.

\begin{pro}[{\cite[Proposition 3.5]{MR4595287}}]\label{PropDualTest}
    The codimension of the image of the multiplication map
    \begin{equation*}
        F_{1110} \otimes V_1^* \rightarrow S^2V_1^* \otimes V_2^* \otimes V_3^*
    \end{equation*}
    is the dimension of the kernel of the skew-symmetrization map 
    \begin{equation*}
        F_{1110}^{\perp} \otimes V_1 \rightarrow \Lambda^2V_1 \otimes V_2 \otimes V_3.
    \end{equation*}
    The analogous property for each map in $(1210), (1120),...,(0112)$-tests holds. In addition, the kernel of the transpose of the map (\ref{eq(1111)test}) is the space
    \begin{equation*}
        (F_{1110}^{\perp} \otimes V_4) \cap (F_{1101}^{\perp} \otimes V_3) \cap (F_{1011}^{\perp} \otimes V_2) \cap (F_{0111}^{\perp} \otimes V_1).
    \end{equation*}
\end{pro}

%-------------------------------------------------------------

\section{Proof of \Cref{ThmMain} and some remarks}\label{Sect:Proof}

As $G_{\det_4} \supset \SL_4$, we take the Lie group $\b$ of upper-triangular matrices with determinant $1$ which is the Borel subgroup of $\SL_4$. Let $\mathfrak{b}$ denote the Borel subalgebra of $\mathfrak{sl}_4$ that is the Lie algebra of $\mathbb{B}$. Since $\mathbb{B}$ is connected, all the $\mathfrak{b}$-fixed subspaces are exactly the $\b$-fixed subspaces.

We decompose the Borel subalgebra as $\mathfrak{b}=\mathfrak{t} \oplus \mathfrak{n}$, where $\mathfrak{t}$ is the algebra of diagonal matrices with trace zero and $\mathfrak{n}$ is the set of strictly upper-triangular matrices. Here, $\mathfrak{t}$ is the Cartan subalgebra of $\mathfrak{sl}_4$. For a given irreducible representation $W$ of $\mathfrak{sl}_4$, we draw the weight diagram of $W$ where the nodes are eigenspaces of the action of the Cartan subalgebra $\mathfrak{t}$ and the arrows indicate where each weight space moves by an action of $\mathfrak{n}$.

Let 
\begin{equation}
	\left\{v_1=\begin{bmatrix} 1 \\ 0 \\ 0 \\ 0\end{bmatrix}, v_2=\begin{bmatrix} 0 \\ 1 \\ 0 \\ 0\end{bmatrix}, v_3=\begin{bmatrix} 0 \\ 0 \\ 1 \\ 0\end{bmatrix}, v_4=\begin{bmatrix} 0 \\ 0 \\ 0 \\ 1\end{bmatrix}
		 \right\}
\end{equation}
be a standard basis of $V:=\c^4$. Let $L_i \in \mathfrak{t}^*$ be defined as
\begin{equation}
	L_i\left( \begin{bmatrix}
		a_1 &  & &\\
		& a_2 & & \\
		 & & a_3 & \\
		 & & & a_4
	\end{bmatrix} \right)=a_i
\end{equation}
for each $i \in {[ 4 ]}$. These are the weights of the standard representation $V$ of $\mathfrak{sl}_4$.

In this section, we prove that there is no ideal $I$ corresponding to border rank $11$ decomposition, by using tests in Section \ref{SubSect:IBorderRankCriterion}. Since $\det_4$ is a consice tensor, we let $I_{1000}=I_{0100}=I_{0010}=I_{0001}=0$ (see Remark \ref{RmkConciseIdeal}). At first, we show that all $(2100),(1200),\cdots,(0012)$ tests must pass in Section \ref{Subsect:Finding $F_{1100}$ and Related Tests}. Skipping $(1110)$, $(1101)$, $(1011)$, $(0111)$-tests, we do $(2110)$, $(1210)$, $(1120)$-tests by finding all $\b$-fixed subspaces $F_{1110}$, and pick a single $\b$-fixed subspace which pass all the tests in Section \ref{Subsect:Finding $F_{1110}$ and Related Tests}. Similarly, we pick three more $\b$-fixed subspaces for  $F_{1101},F_{1011}$ and $F_{0111}$ to do corresponding tests in the same subsection. Finally, using the candidates which passed these tests, we do $(1111)$-test in Section \ref{$(1111)$-test}. This $(1111)$-test gives us the conclusion that $\underline{\rk}(\det_4) > 11$.
%In Section \ref{Subsect:Some remarks}, we give some remarks on the fixed ideal theorem and a question.
In Section \ref{Subsect:Some remarks}, we give the final remarks.

\subsection{$F_{1100}$ and Related Tests}\label{Subsect:Finding $F_{1100}$ and Related Tests}

Let  
\begin{equation}
	F_{1100} \subset \det_4(V_3^* \otimes V_4^*)^{\perp} \subset V_1^* \otimes V_2^*
\end{equation}
be a $\b$-fixed subspace of codimension $11$ in $V_1^* \otimes V_2^*$, i.e., of dimension $16-11=5$. Then the map
\begin{equation*}
	F_{1100} \otimes V_1^* \rightarrow S^2V_1^* \otimes V_2^*
\end{equation*}
has rank at most $20$, i.e., has image of codimension at least $\dim(S^2V_1^* \otimes V_2^*)-20=20$. Since the image has codimension greater than $11$, all such $F_{1100}$ automatically pass $(2100)$-test. Similarly, $F_{1100}$ must pass $(1200)$-test. Moreover, all the candidates $F_{1010},F_{1001},F_{0110},F_{0101},F_{0011}$ must pass $(2010),(1020)$-tests, $(2001),(1002)$-tests, $(0210),(0120)$-tests, $(0201),(0102)$-tests, $(0021),(0012)$-tests, respectively.

%The following is the weight diagram of $S^2V$. Here, we define $u(i,j)=(v_i\otimes v_j)\cdot c_{(2)}$ where $c_{(2)}$ is the symmetrizer corresponding to the partition $(2)$ of 2. Note that it has the weight $L_i+L_j$.
%\[\begin{tikzcd}
%	&u(1,1) & \\
%	&u(1,2)\arrow[u] & \\
%	u(1,3)\arrow[ur]& & u(2,2)\arrow[ul]\\
%	u(1,4)\arrow[u]& & u(2,3)\arrow[ull]\arrow[u]\\
%	u(2,4)\arrow[u]\arrow[urr]&&u(3,3)\arrow[u]\\
%	&u(3,4)\arrow[ul]\arrow[ur]&\\
%	&u(4,4)\arrow[u]&
%\end{tikzcd}\]

\subsection{Finding $F_{1110}$ and Related Tests}\label{Subsect:Finding $F_{1110}$ and Related Tests}

Here, we find all $\b$-fixed (i.e., $\mathfrak{b}$-fixed) subspaces 
\begin{equation}
	F_{1110} \subset \det_4(V_4^*)^{\perp} \subset V_1^* \otimes V_2^* \otimes V_3^*
\end{equation}
of codimension $11$ in $V_1^* \otimes V_2^* \otimes V_3^*$, i.e., of dimension $64-11=53$. We will find $F_{1110}$ by finding 
\begin{equation*}
	E_{1110}:=F_{1110}^{\perp} \subset V_1 \otimes V_2 \otimes V_3
\end{equation*}
which is $\mathbb{B}$-fixed subspace of dimension $11$. Recall that $V_1 = V_2 = V_3 = V_4 = \c^4 =: V$.  Note that
\begin{equation*}
	\Lambda^3V=\det_4(V_4^*) \subset F_{1110}^{\perp}=E_{1110}.
\end{equation*} 
and $\dim \Lambda^3V=4$.
Define Young tableaux as
\begin{equation*}
	\lambda_1=
	{\scriptsize\begin{ytableau}
			1 & 2  & 3
	\end{ytableau}}
	~~
	\text{, }
	~~
	\lambda_2=
	{\scriptsize\begin{ytableau}
			1 & 2  \\
			3
	\end{ytableau}}
	~~
    \text{, }
    ~~
	\lambda_3=
	{\scriptsize\begin{ytableau}
			1 & 3  \\
			2
	\end{ytableau}}
	~~
	\text{, and}
	~~
	\lambda_4=
	{\scriptsize\begin{ytableau}
	1  \\
	2 \\
    3
	\end{ytableau}}
	~~
	.
\end{equation*}
Then
\begin{equation*}
	V^{\otimes 3} \cong S_{\lambda_1}V \oplus S_{\lambda_2}V \oplus S_{\lambda_3}V \oplus S_{\lambda_4}V
\end{equation*}
as a representation of $\operatorname{SL}_4$ where $S_{\lambda_1}V=S^3V$ and $S_{\lambda_4}V=\Lambda^3V$.
Hence we need to find $7$-dimensional $\mathbb{B}$-fixed subspaces
\begin{equation*}
    E'_{1110}\subset S_{\lambda_1}V\oplus S_{\lambda_2}V\oplus S_{\lambda_3}V
\end{equation*}
to get $E_{1110}=E'_{1110}\oplus \Lambda^3V$.
The Young symmetrizers are
\begin{align*}
    c_{\lambda_1} &= (1)+(1\ 2)+(1\ 3)+(2\ 3)+(1\ 2\ 3)+(1\ 3\ 2)\\
    c_{\lambda_2} &= ((1)+(1\ 2))((1)-(1\ 3))\\
    c_{\lambda_3} &= ((1)+(1\ 3))((1)-(1\ 2))\\
    c_{\lambda_4} &= (1)-(1\ 2)-(1\ 3)-(2\ 3)+(1\ 2\ 3)+(1\ 3\ 2).
\end{align*}
We denote $$u_l(i,j,k):=(v_i\otimes v_j\otimes v_k)\cdot c_{\lambda_l}$$
which has the weight $L_i+L_j+L_k$.
Then we get the following weight diagram of $S_{\lambda_1}V\oplus S_{\lambda_2}V\oplus S_{\lambda_3}V$.

\begin{figure}[H]
\adjustbox{max width=\textwidth}{
\begin{tikzcd}[row sep=small]
	&&\mlnode{$u_1(1,1,1)$} && \\
	&&\mlnode{$u_1(1,1,2)$\\$u_2(1,1,2)$\\$u_3(1,1,2)$}\arrow[u] && \\
	&\mlnode{$u_1(1,1,3)$\\$u_2(1,1,3)$\\$u_3(1,1,3)$}\arrow[ur]& &\mlnode{$u_1(1,2,2)$\\$u_2(2,2,1)$\\$u_3(2,2,1)$}\arrow[ul]&\\
	\mlnode{$u_1(1,1,4)$\\$u_2(1,1,4)$\\$u_3(1,1,4)$}\arrow[ur]& & \mlnode{$u_1(1,2,3)$\\$u_2(1,2,3)$\\$u_2(1,2,3)+u_2(1,3,2)$\\$u_3(1,3,2)$\\$u_3(1,2,3)+u_3(1,3,2)$}\arrow[ur]\arrow[ul]&&\mlnode{$u_1(2,2,2)$}\arrow[ul]\\
	\mlnode{$u_1(1,2,4)$\\$u_2(1,2,4)$\\$u_2(1,2,4)+u_2(1,4,2)$\\$u_3(1,4,2)$\\$u_3(1,2,4)+u_3(1,4,2)$}\arrow[u]\arrow[urr]& & \mlnode{$u_1(1,3,3)$\\$u_2(3,3,1)$\\$u_3(3,3,1)$}\arrow[u]&&\mlnode{$u_1(2,2,3)$\\$u_2(2,2,3)$\\$u_3(2,2,3)$}\arrow[ull]\arrow[u]\\
	\mlnode{$u_1(1,3,4)$\\$u_2(1,3,4)$\\$u_2(1,3,4)+u_2(1,4,3)$\\$u_3(1,4,3)$\\$u_3(1,3,4)+u_3(1,4,3)$}\arrow[u]\arrow[urr]&&\mlnode{$u_1(2,2,4)$\\$u_2(2,2,4)$\\$u_3(2,2,4)$}\arrow[ull]\arrow[urr]&&\mlnode{$u_1(2,3,3)$\\$u_2(3,3,2)$\\$u_3(3,3,2)$}\arrow[ull]\arrow[u]\\
	\mlnode{$u_1(1,4,4)$\\$u_2(4,4,1)$\\$u_3(4,4,1)$}\arrow[u]&&\mlnode{$u_1(2,3,4)$\\$u_2(2,3,4)$\\$u_2(2,3,4)+u_2(2,4,3)$\\$u_3(2,4,3)$\\$u_3(2,3,4)+u_3(2,4,3)$}\arrow[ull]\arrow[u]\arrow[urr]&&\mlnode{$u_1(3,3,3)$}\arrow[u]\\
	&\mlnode{$u_1(2,4,4)$\\$u_2(4,4,2)$\\$u_3(4,4,2)$}\arrow[ur]\arrow[ul]& &\mlnode{$u_1(3,3,4)$\\$u_2(3,3,4)$\\$u_3(3,3,4)$}\arrow[ul]\arrow[ur]&\\
	&&\mlnode{$u_1(3,4,4)$\\$u_2(4,4,3)$\\$u_3(4,4,3)$}\arrow[ur]\arrow[ul]&&\\
	&&\mlnode{$u_1(4,4,4)$}\arrow[u]&&\\
\end{tikzcd}
}
\caption{The weight diagram of $S_{\lambda_1}V\oplus S_{\lambda_2}V\oplus S_{\lambda_3}V$}
\end{figure}

We denote these 20 weight spaces by $W_1,W_2,\cdots,W_{20}$.
Let $d_i$ be the dimension of $E
'_{1110}\cap W_i$ for each $i$.
Then we have $d_1+d_2+\cdots+d_{20}=7$.
Also, whenever we have the arrows $W_i\to W_{j_1}$, $W_i\to W_{j_2}$, $\cdots$, $W_i\to W_{j_k}$, we get the condition
\[
    d_{j_1}+d_{j_2}+\cdots+d_{j_k}+\dim\ker\left(W_i\to W_{j_1}\oplus W_{j_2}\oplus\cdots\oplus W_{j_k}\right)\geq d_{i}.
\]
There are 228 sequences of $(d_1,d_2,\cdots,d_{20})$ satisfying these conditions.
For each $(d_1,d_2,\cdots,d_{20})$, we consider the subspace corresponding to an element in the product of Grassmannians 
\[
G(d_1,w_1)\times G(d_2,w_2)\times \cdots\times G(d_{20},w_{20})
\]
where $w_i:=\dim W_i$.
Then we represent an element in $G(d_i,w_i)$ as a $d_i\times w_i$ matrix and decompose $G(d_i,w_i)$ into $\binom{w_i}{d_i}$ subsets as follows:
\begin{enumerate}
    \item for a given $\{n_1,n_2,\cdots,n_{d_i}\}\subset \{1,2,\cdots,w_i\}$ with $n_1< n_2< \cdots < n_{d_i}$, the $(j, n_j)$-th entry is 1 and all entries lying to its left in the same row or lying in the same column are zero;
    \item\label{item:variable} the other entries are represented by distinct variables.
\end{enumerate}
Then these subsets are pairwise disjoint and their union is $G(d_i,w_i)$.

For example, every point of $G(2,4)$ is the row space of exactly one of the following 6 matrices.

\begin{align*}
&{\begin{bmatrix}
    1&0&x_0&x_1\\
    0&1&x_2&x_3
\end{bmatrix}}
&&{\begin{bmatrix}
    1&x_0&0&x_1\\
    0&0&1&x_2
\end{bmatrix}}
&&{\begin{bmatrix}
    1&x_0&x_1&0\\
    0&0&0&1
\end{bmatrix}}
\\
&{\begin{bmatrix}
    0&1&0&x_0\\
    0&0&1&x_1
\end{bmatrix}}
&&{\begin{bmatrix}
    0&1&x_0&0\\
    0&0&0&1
\end{bmatrix}}
&&{\begin{bmatrix}
    0&0&1&0\\
    0&0&0&1
\end{bmatrix}}
\end{align*}

This requires us to consider $\binom{w_1}{d_1}\times \binom{w_2}{d_2}\times \cdots \times \binom{w_{20}}{d_{20}}$ cases for each $(d_1,d_2,\cdots,d_{20})$.
For each case, we compute the locus of the parameters in \eqref{item:variable} such that the assembled subspace is $\b$-fixed.
Here, we get 322 parametrized families of $\b$-fixed subspaces, whose members are the candidates for $E'_{1110}$. 
Using Proposition \ref{PropDualTest}, we checked whether each $E_{1110}=E'_{1110}\oplus \Lambda^3 V$ passes (2110),(1210),(1120)-tests.
There are 128 parametrized families for the candidates of $E_{1110}$ that pass the (2110)-test.
For the survived ones, there are 11 parametrized families surviving the (1210)-test, and there is only one candidate surviving the (1120)-test simultaneously.
This unique candidate is
	\begin{align*}
		E_{1110}=&\langle u_1(1,1,1)\rangle\oplus\langle u_1(1,1,2),u_2(1,1,2),u_3(1,1,2)\rangle\\
        &\oplus \langle u_1(1,1,3),u_2(1,1,3),u_3(1,1,3)\rangle\oplus \Lambda^3V.
	\end{align*}

By the symmetry, we determine the unique candidate for each of $E_{1101}, E_{1011}$, and $E_{0111}$. 

\subsection{$(1111)$-test}\label{$(1111)$-test}

For this $4$-tuple $\{E_{1110},E_{1101},E_{1011},E_{0111}\}$, the space
\begin{equation*}
	(E_{1110} \otimes V_4) \cap (E_{1101} \otimes V_3)\cap(E_{1011} \otimes V_2)\cap(E_{0111} \otimes V_1)
\end{equation*}
has dimension 10.
This space is the kernel of the transpose of the map
\begin{equation}\label{eqn:1111}
	(F_{1110} \otimes V_4^*) \oplus (F_{1101} \otimes V_3^*)\oplus(F_{1011} \otimes V_2^*)\oplus(F_{0111} \otimes V_1^*) \rightarrow V_1^* \otimes V_2^* \otimes V_3^* \otimes V_4^*
\end{equation}
by Proposition \ref{PropDualTest}.
Hence the image of the map \eqref{eqn:1111} has codimension 10.
This proves $\brk(\det_4)>11$.

Consequently, we conclude $\brk(\det_4)=12$ over $\c$, and so over any subfield of $\c$.

%In \cite[Remark 4.4]{MR4332674}, Buczy\'nska--Buczy\'nski expected the fixed ideal theorem (\Cref{thmFixedIdeal}) does not hold if one replace the border rank to the rank. In this section, we show there indeed exists a tensor of rank $r$ that does not have a decomposition into $r$ rank one tensors that corresponds to a $\b$-fixed ideal $I$ where $\b$ is a connected solvable group contained in the symmetric group of the tensor. Here, the ideal corresponding to the decomposition is defined similarly as before.

%\subsection{Some remarks on the fixed ideal theorem}\label{Subsect:Some remarks}
\subsection{Final remarks}\label{Subsect:Some remarks}
Note that $\rk(\det_n)=\brk(\det_n)$ for $n\leq 4$. So one may ask the following question:
\begin{question}
	$\mathbf{R}(\det_n)=\underline{\mathbf{R}}(\det_n)$ for all $n \in \mathbb{N}$?
\end{question}
Answering this question seems very challenging, since both the tensor rank and the border rank are far from being determined for $n\geq 5$.
To the authors' knowledge, the best known bounds of $\rk(\det_5)$ and $\brk(\det_5)$ are
\[
	27\leq \brk(\det_5)\leq \rk(\det_5)\leq 52
\]
where the lower bound is obtained from \cite[Theorem 7]{RKF} and the upper bound is obtained from \cite[Corollary 3]{MR4822414}.

\begin{rem}\label{rem: det2}
	In \cite[Remark 4.4]{MR4332674}, Buczy\'nska--Buczy\'nski expected the tensor rank version of the fixed ideal theorem (\Cref{thmFixedIdeal}) does not hold.
	In other words, they expected that there exists a tensor $T\in V_1\otimes V_2\otimes \cdots \otimes V_d$ that does not yield any decomposition $T=T_1+\cdots +T_r$ into rank one tensors such that the multigraded ideal $I\subset \sym(V_1^*)\otimes \cdots\otimes \sym(V_d^*)$ corresponding to $T_1,\cdots,T_r$ satisfies all of the following:
	\begin{itemize}
		\item [(\romannumeral1)] $I \subset \operatorname{Ann}(T)$
		\item [(\romannumeral2)] For all $(i_1,\cdots,i_d)$ such that $r \leq \dim (S^{i_1}V_1 ^*\otimes \cdots \otimes S^{i_d}V_d^*)$, we have $$\codim (I_{i_1,\cdots,i_d}) = r.$$
		\item [(\romannumeral3)] The multiplication map
		{\small
		\begin{equation*}
			(I_{(i_1-1),\cdots,i_d} \otimes V_1^*) \oplus \cdots \oplus (I_{i_1,\cdots,(i_d-1)} \otimes V_d^*)
			\rightarrow S^{i_1}V_1^*\otimes \cdots \otimes S^{i_d}V_d^*
		\end{equation*}
		}
		has its image in $I_{i_1,\cdots,i_d}$. 
		\item [(\romannumeral4)] Each $I_{i_1,\cdots,i_d}$ is $\mathbb{B}$-fixed.
	\end{itemize}
	In this remark, we verify their expectation.
	Let $V=V_1=V_2=\c^2$ and $T=\det_2=e_1 \otimes e_2-e_2 \otimes e_1\in V_1\otimes V_2$, where $\{e_1,e_2\}$ is the standard basis of $\mathbb{C}^2$. As $\SL_2\subset G_T$, we let $\b$ be the Borel subgroup of $\SL_2$. Note that $\brk(\det_2)=2$.
	Suppose that there is a decomposition of $\det_2$ into two rank one tensors whose corresponding ideal $I\subset \sym(V_1 ^*)\otimes \sym(V_2^*)$ satisfies all the conditions (i)-(iv).
		Let $F_{11}$ be the candidates of $I_{11}$.
		Recall that $V^*\otimes V^*\cong S^2 V^*\oplus \Lambda^2 V^*$ as $\SL_2$-representations. The weight diagrams of $S^2 V^*$ and $\Lambda^2 V^*$ are as follows.
\begin{figure}[H]
\[\begin{tikzcd}
		x^2\arrow[d]&\\
		xy\arrow[d]&\hspace{30pt} x\wedge y\\
		y^2&
	\end{tikzcd}
    \]
    \caption{Weight diagrams of $S^2V^*$ and $\Lambda^2V^*$}
\end{figure}
    \noindent
    where $\{x,y\}$ is the dual basis of $\{e_1,e_2\}$. Hence, there is a one-dimensional family of $\mathbb{B}$-fixed subspaces of codimension $2$ in $V^* \otimes V^*$, which is
    \begin{equation*}
        \{\langle y \otimes y, (p+q)x \otimes y + (p-q) y \otimes x \rangle~|~[p:q] \in \mathbb{P}^1\}.
    \end{equation*}
    Since 
    \begin{equation*}
        F_{11} \subset \operatorname{Ann}(\operatorname{det}_2)_{11}=\langle x \otimes x, y \otimes y, x \otimes y + y \otimes x \rangle,
    \end{equation*}
    then there exists a unique candidate
    \begin{equation*}
        F_{11}=\langle y \otimes y, x \otimes y + y \otimes x \rangle.
    \end{equation*}
	One can easily see that the vanishing locus of this candidate is a singleton $\{([1:0],[1:0])\}$.
	In other words, no such an ideal $I$ corresponds to a rank decomposition of $\det_2$.
\end{rem}

\bibliography{Border}
\vspace{0.5cm}

\end{document}